\documentclass[12pt]{amsart}
\usepackage{amssymb,latexsym}
\newcommand{\ben}{\begin{enumerate}}
\newcommand{\een}{\end{enumerate}}
\newcommand{\ble}{\begin{lem}}
\newcommand{\ele}{\end{lem}}
\newcommand{\bth}{\begin{thm}}
\renewcommand{\eth}{\end{thm}}
\newcommand{\bpr}{\begin{prop}}
\newcommand{\epr}{\end{prop}}
\newcommand{\bco}{\begin{cor}}
\newcommand{\eco}{\end{cor}}
\newcommand{\bcon}{\begin{conj}}
\newcommand{\econ}{\end{conj}}
\newcommand{\bde}{\begin{defn}}
\newcommand{\ede}{\end{defn}}
\newcommand{\bex}{\begin{exa}}
\newcommand{\eex}{\end{exa}}
\newcommand{\barr}{\begin{array}}
\newcommand{\earr}{\end{array}}
\newcommand{\btab}{\begin{tabular}}
\newcommand{\etab}{\end{tabular}}
\newcommand{\beq}{\begin{equation}}
\newcommand{\eeq}{\end{equation}}
\newcommand{\bea}{\begin{eqnarray*}}
\newcommand{\eea}{\end{eqnarray*}}
\newcommand{\bce}{\begin{center}}
\newcommand{\ece}{\end{center}}
\newcommand{\bpi}{\begin{picture}}
\newcommand{\epi}{\end{picture}}
\newcommand{\bfi}{\begin{figure} \begin{center}}
\newcommand{\efi}{\end{center} \end{figure}}

\newcommand{\bsl}{\begin{slide}{}}
\newcommand{\esl}{\end{slide}}

\newcommand{\bib}{thebibliography}

\newcommand{\ol}{\overline}

\newcommand{\hso}[1]{\hspace{-1pt}}

\newcommand{\qmq}[1]{\quad\mbox{#1}\quad}

\newcommand{\ptn}{\vdash}

\newcommand{\case}[4]{\left\{\barr{ll}#1&\mbox{#2}\\#3&\mbox{#4}\earr\right.}

\def\<{\langle}
\def\>{\rangle}

\newcommand{\ree}[1]{(\ref{#1})}

\newcommand{\ra}{\rightarrow}

\newcommand{\lra}{\longrightarrow}

\newcommand{\llra}{\longleftrightarrow}

\newcommand{\ep}{\epsilon}

\newcommand{\ka}{\kappa}
\newcommand{\la}{\lambda}

\newcommand{\La}{\Lambda}

\newcommand{\bx}{{\bf x}}

\newcommand{\bbP}{{\mathbb P}}

\newcommand{\Sb}{\ol{S}}

\newcommand{\Tb}{\ol{T}}

\newcommand{\sh}{\mathop{\rm sh}\nolimits}

\newcommand{\aim}{Adv.\ in Math.\/}

\newcommand{\jcta}{J. Combin.\ Theory Ser. A\/}

\newcommand{\oup}{Oxford University Press}

\newcommand{\Gaa}{\put(0,0){\circle*{3}}}
\newcommand{\Gba}{\put(10,0){\circle*{3}}}
\newcommand{\Gca}{\put(20,0){\circle*{3}}}

\newcommand{\Gab}{\put(0,10){\circle*{3}}}
\newcommand{\Gbb}{\put(10,10){\circle*{3}}}
\newcommand{\Gcb}{\put(20,10){\circle*{3}}}

\newcommand{\Gac}{\put(0,20){\circle*{3}}}
\newcommand{\Gbc}{\put(10,20){\circle*{3}}}
\newcommand{\Gcc}{\put(20,20){\circle*{3}}}

\newcommand{\Gad}{\put(0,30){\circle*{3}}}
\newcommand{\Gbd}{\put(10,30){\circle*{3}}}
\newcommand{\Gcd}{\put(20,30){\circle*{3}}}
\newcommand{\Gdd}{\put(30,30){\circle*{3}}}

\newcommand{\Gae}{\put(0,40){\circle*{3}}}
\newcommand{\Gbe}{\put(10,40){\circle*{3}}}
\newcommand{\Gce}{\put(20,40){\circle*{3}}}

\newcommand{\GbaV}[1]{\put(10,0){\circle*{#1}}}
\newcommand{\GcaV}[1]{\put(20,0){\circle*{#1}}}

\newcommand{\GabV}[1]{\put(0,10){\circle*{#1}}}
\newcommand{\GbbV}[1]{\put(10,10){\circle*{#1}}}
\newcommand{\GcbV}[1]{\put(20,10){\circle*{#1}}}

\newcommand{\GacV}[1]{\put(0,20){\circle*{#1}}}

\newcommand{\GccV}[1]{\put(20,20){\circle*{#1}}}

\newcommand{\GadV}[1]{\put(0,30){\circle*{#1}}}
\newcommand{\GbdV}[1]{\put(10,30){\circle*{#1}}}
\newcommand{\GcdV}[1]{\put(20,30){\circle*{#1}}}
\newcommand{\GddV}[1]{\put(30,30){\circle*{#1}}}

\newcommand{\Gaaba}{\put(0,0){\line(1,0){10}}}

\newcommand{\Gaaab}{\put(0,0){\line(0,1){10}}}

\newcommand{\Gaaac}{\put(0,0){\line(0,1){20}}}

\newcommand{\Gbaca}{\put(10,0){\line(1,0){10}}}

\newcommand{\Gbabb}{\put(10,0){\line(0,1){10}}}

\newcommand{\Gcacb}{\put(20,0){\line(0,1){10}}}

\newcommand{\Gcaac}{\put(20,0){\line(-1,1){20}}}

\newcommand{\Gcacc}{\put(20,0){\line(0,1){20}}}

\newcommand{\Gabbb}{\put(0,10){\line(1,0){10}}}

\newcommand{\Gabac}{\put(0,10){\line(0,1){10}}}

\newcommand{\Gbbcb}{\put(10,10){\line(1,0){10}}}

\newcommand{\Gbbbc}{\put(10,10){\line(0,1){10}}}

\newcommand{\Gcbcc}{\put(20,10){\line(0,1){10}}}

\newcommand{\Gacbc}{\put(0,20){\line(1,0){10}}}

\newcommand{\Gacad}{\put(0,20){\line(0,1){10}}}

\newcommand{\Gbcbd}{\put(10,20){\line(0,1){10}}}

\newcommand{\Gccdc}{\put(20,20){\line(1,0){10}}}

\newcommand{\Gadbd}{\put(0,30){\line(1,0){10}}}

\newcommand{\Gbdcd}{\put(10,30){\line(1,0){10}}}

\newcommand{\Gbdbe}{\put(10,30){\line(0,1){10}}}

\newcommand{\Gcddd}{\put(20,30){\line(1,0){10}}}

\newcommand{\Gdded}{\put(30,30){\line(1,0){10}}}

\newcommand{\Gbece}{\put(10,40){\line(1,0){10}}}

\newtheorem{thm}{Theorem}[section]
\newtheorem{prop}[thm]{Proposition}
\newtheorem{cor}[thm]{Corollary}
\newtheorem{lem}[thm]{Lemma}
\newtheorem{conj}[thm]{Conjecture}
\newtheorem{exa}[thm]{Example}

\newcommand{\tpo}{{\rm (}{\bf 3}{\rm +}{\bf 1}{\rm )}}
\newcommand{\Sd}{\dot{S}}

\newcommand{\sqr}{\rule{8pt}{8pt}}

\newcommand{\Par}{\mathop{\rm Par}\nolimits}
\newcommand{\h}{\mbox{---}}

\begin{document}
\title{A sign-reversing involution for rooted special rim-hook tableaux
}
\author{Bruce E. Sagan}
\address{Department of Mathematics, Michigan State University,
East Lansing, MI 48824-1027, USA}
\email{sagan@math.msu.edu}

\author{Jaejin Lee}
\address{Department of Mathematics, Hallym University,
Chunchon, KOREA 200-702}
\email{jjlee@hallym.ac.kr}

\date{\today}

	\begin{abstract}
E\~gecio\~glu and Remmel~\cite{er:cii} gave an interpretation for the
entries of the inverse Kostka matrix $K^{-1}$ in terms of special rim-hook
tableaux.  They were able to use this interpretation to give a
combinatorial proof that $KK^{-1}=I$ but were unable to do the same
for the equation $K^{-1}K=I$.  We define a sign-reversing involution
on rooted special rim-hook tableaux which can be used to prove that 
the last column of this second product is correct.  In addition,
following a suggestion of Chow~\cite{cho:nci}  we combine our
involution with a result of Gasharov~\cite{gas:igt} to give
a combinatorial proof  of a  special case of the \tpo-free
Conjecture of Stanley and Stembridge~\cite{ss:ijt}. 
	\end{abstract}

\maketitle

\section{Introduction} 
\label{i}

We first recall some definitions from the theory of Young tableaux.
Further exposition can be found in the texts of Fulton~\cite{ful:yt},
Macdonald~\cite{mac:sfh}, Sagan~\cite{sag:sg}, and
Stanley~\cite{sta:ec2}.  

Let $\la=(\la_1,\la_2,\ldots,\la_l)$ be a partition of the nonnegative
integer $n$, denoted $\la\ptn n$, so $\la$ is a weakly decreasing
sequence of positive integers summing to $n$.  We will also let $\la$
stand for the Ferrers diagram of $\la$ written in English notation
with $\la_i$ nodes or cells in the $i$th row from the top.  In
addition, we will 
use the notation $\la=(1^{m_1},2^{m_2},\ldots,n^{m_n})$ which means
that the integer $j$ appears $m_j$ times in $\la$.  The set of all
partitions of all $n\ge0$ will be denoted $\Par$. 

Given $\la,\mu\ptn n$,
the corresponding {\it Kostka number\/} $K_{\la,\mu}$ is the number of
semistandard Young tableaux $T$ of shape $\sh(T)=\la$ and content
$c(T)=\mu$, i.e., the number of fillings of the diagram $\la$ with
integers such that 
rows weakly increase, columns strictly increase, and for all $k$ the integer $k$
occurs $\mu_k$ times.  For fixed $n$, we collect these numbers into
the {\it Kostka matrix\/} $K=(K_{\la,\mu})$.  We will use the
reverse lexicographic order on  partitions so that $K$ becomes
upper unitriangular.

E\~gecio\~glu and Remmel~\cite{er:cii} gave a combinatorial
interpretation for the 
entries of the inverse Kostka matrix $K^{-1}$ as follows.
A {\it rim hook}, $H$, is a skew diagram (the set-theoretic difference
of two ordinary diagrams) which is connected and contains no
$2\times2$ square of cells.  The {\it size of $H$\/} is the number
of cells it contains.  A {\it rim-hook tableau}  of shape
$\la$ is a partition of the diagram of $\la$ into rim hooks.  This
tableau $S$ is {\it special\/} if each of the rim hooks contains a cell
from the first column of $\la$.  The {\it type\/} of $S$ is
$t(S)=(1^{m_1},2^{m_2},\ldots,n^{m_n})$ where $m_k$ is the number of
rim hooks in $S$ of size $k$.
Using nodes for the Ferrers diagram and connecting them if they are
adjacent in
the same rim hook, the following diagram illustrates a special rim-hook
tableau $S$ with $\sh(S)=(1^2,2^2,3)$ and  $t(S)=(1,4^2)$.
\medskip
$$
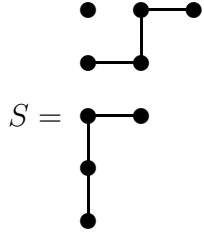
\thicklines
\setlength{\unitlength}{2pt}
\bpi(20,40)(0,0)
\put(-10,20){\makebox(0,0){$S=$}}
\Gaa \Gab \Gac \Gad \Gae
\Gbc \Gbd \Gbe
\Gce
\Gaaab \Gabac \Gacbc
\Gadbd \Gbdbe \Gbece
\epi
$$ 
\medskip

The {\it leg length of rim hook $H$}, $l(H)$, is the number of
vertical edges in $H$ when viewed as in the diagram above.
We now define the {\it sign\/} of a rim
hook $H$ and of a rim-hook tableau $S$ to be
$$
\ep(H)=(-1)^{l(H)}\qmq{and} \ep(S)=\prod_{H\in S} \ep(H),
$$
respectively.  The previous rim-hook tableau has sign
$$
\ep(S)=(-1)^{0}\cdot(-1)^{1}\cdot(-1)^{2}=-1.
$$
We can now state E\~gecio\~glu and Remmel's interpretation.
\bth[E\~gecio\~glu and Remmel]
\label{K-1}
The entries of the inverse Kostka matrix are given by
$$
K^{-1}_{\mu,\la}=\sum_{S} \ep(S)
$$
where the sum is over all special rim hook tableaux $S$ with $\sh(S)=\la$
and $t(S)=\mu$.\qed
\eth

In~\cite{er:cii} they used this theorem to give a combinatorial proof
of the fact that $KK^{-1}=I$ but were not able to do the same thing
for the identity $K^{-1}K=I$.  In the next section, we will give an
algorithmic, sign-reversing involution which will prove that the last
column of the second product is correct.  
Parts of this procedure are reminiscent of the lattice path involution
of Lindstr\"om~\cite{lin:vri} and Gessel-Viennot~\cite{gv:bdp,gv:dpp}
as well as the rim-hook Robinson-Schensted algorithm of
White~\cite{whi:bpo} and Stanton-White~\cite{sw:scr}. 
In section~\ref{tfc} we follow a suggestion
of Chow~\cite{cho:nci} and combine our involution with a result of
Gasharov~\cite{gas:igt} to prove 
a special case of the \tpo-free Conjecture of Stanley and
Stembridge~\cite{ss:ijt}.  Finally, we end with a discussion of
further work which needs to be done.

\section{The basic involution}
\label{bi}

First note that by Theorem~\ref{K-1}, $K^{-1}_{(1^n),\la}$ is just
the number of special rim-hook tableaux of shape $\la$ where all hooks
have size one (since such tableaux have sign $+1$).  But
since they must also contain a cell of the first column, such a
tableau exists precisely when $\la=(1^n)$ and in that case the tableau
is unique.  Since the definition of the Kostka
matrix immediately gives $K_{(1^n),(1^n)}=1$, we have that the inner
product of the last row of $K^{-1}$ and the last column of $K$ is $1$
as desired.  So for the rest of this section we will assume that the
type of our special rim-hook tableau is $\mu\neq (1^n)$.

We wish to show that the inner product of the $\mu$th row of $K^{-1}$ with
column $(1^n)$ of $K$ is $0$.  It follows from E\~gecio\~glu and
Remmel's Theorem that we need to prove
\beq
\label{ST}
\sum_{(S,T)} \ep(S)=0
\eeq
the sum being over all pairs $(S,T)$ where $S$ is a special rim-hook
tableau with $t(S)=\mu$, and $T$ is a standard Young tableau (that
is,  $c(T)=(1^n)$) of the same shape as $S$.  We will prove this
identity by exhibiting a sign-reversing involution $I$ on such pairs.

Suppose first that the cell of $n$ in $T$ corresponds to a hook of
size one in $S$.  Then since $S$ is special, this cell is at the end
of the first column.  In this case, remove that cell from both $S$
and $T$ to form $\Sb$ and $\Tb$ respectively.  Now we can assume, by
induction on $n$, that $I(\Sb,\Tb)=(\Sb',\Tb')$ has been defined.  So
let $I(S,T)=(S',T')$ where $S'$ is $\Sb'$ with a hook of size $1$
added to the end of the first column  and  $T'$ is $\Tb'$ with a cell
labeled $n$ added  to the end of the first column.  Clearly this will
result in a sign-reversing involution as long as this was true for
pairs with $n-1$ cells.  So for the rest of this section we will
also assume that the cell containing $n$ in $T$ corresponds to a cell in a
hook of at least two cells in $S$.

To describe $I$ under these assumptions, we will need names for the
different parts of a 
rim hook $H$.  As usual, let $(i,j)$ denote the cell of a shape $\la$
in row $i$ and column $j$.  An {\it internal corner\/} of $H$ is
$(i,j)\in H$ such that $(i+1,j),(i,j+1)\in H$.   
Dually, an {\it external corner\/} of $H$ is $(i,j)\in H$ such
that $(i-1,j),(i,j-1)\in H$.   
The {\it head\/} of $H$ is the
$(i,j)\in H$ with smallest $i$ and largest $j$.  
Similarly, the {\it tail\/} of $H$ is the
$(i,j)\in H$ with largest $i$ and $j=1$.  
In the previous example, the upper hook of size 4 has internal corner
$(1,2)$, external corner $(2,2)$, head $(1,3)$, and tail $(2,1)$.
The {\it permissible cells\/} of $H$ are precisely those which are
either an internal corner, external corner, head, or tail. 
 
Now define a {\it rooted Ferrers diagram\/} to 
be a Ferrers diagram where one of the nodes has been marked.  We will
indicate this in our figures by making the distinguished
node a square.   Any tableau built out of a Ferrers
diagram can be rooted in an analogous way.  
Any invariants of the original tableau will be carried over to the 
rooted version without change.
Now associate with any
pair $(S,T)$ in the sum~\ree{ST} a rooted special rim-hook tableau
$\Sd$ by rooting $S$ at the node where the entry $n$ occurs in $T$.
We will define a sign-reversing involution $\iota$ on the set of
rooted special rim-hook tableaux of given type which are obtainable in
this way. 
In addition, $\iota$ will have the property that if $\iota(\Sd)=\Sd'$ and 
$\Sd,\Sd'$ have roots $r,r'$ respectively, then
\beq
\label{sh}
\sh(\Sd)-r=\sh(\Sd')-r'
\eeq
where the minus sign represents set-theoretic difference of diagrams.  Our full
involution $I(S,T)=(S',T')$ will then be the composition 
$$
(S,T)\lra \Sd \stackrel{\iota}{\lra} \Sd'\lra (S',T')
$$
where $S'$ is obtained from $\Sd'$ by forgetting about the root and
$T'$ is obtained by replacing the root of $\Sd'$ by $n$ and leaving
the numbers $1,2,\ldots,n-1$ in the same positions as they were in
$T$.  Note that~\ree{sh} guarantees that $T'$ is well defined.
Furthermore, it is clear from construction that $I$ will be a
sign-reversing involution because $\iota$ is.  Even though $\iota$ has not
been fully defined, we can give an example of the rest of the algorithm as
follows.
\medskip
$$
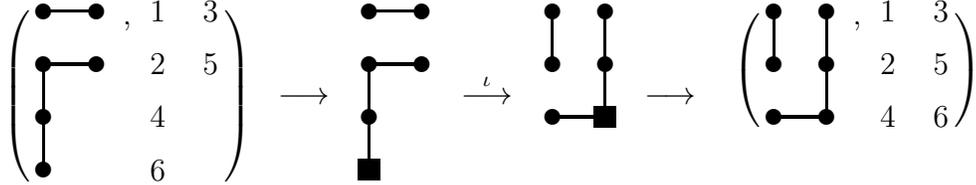
\thicklines
\setlength{\unitlength}{2pt}
\barr{ccccccc}
\left(\barr{ccc}
\bpi(10,30)(0,0)
\Gaa \Gab \Gac \Gad \Gbc \Gbd
\Gaaab \Gabac \Gacbc 
\Gadbd
\epi
&\raisebox{55pt}{,}&
\bpi(10,30)(0,0)
\put(0,30){\makebox(0,0){$1$}}
\put(0,20){\makebox(0,0){$2$}}
\put(10,30){\makebox(0,0){$3$}}
\put(0,10){\makebox(0,0){$4$}}
\put(10,20){\makebox(0,0){$5$}}
\put(0,0){\makebox(0,0){$6$}}
\epi
\earr\right)
&\lra&
\barr{c}
\bpi(10,30)(0,0)
\put(0,0){\makebox(0,0){$\sqr$}}
\Gab \Gac \Gad \Gbc \Gbd
\Gaaab \Gabac \Gacbc 
\Gadbd
\epi
\earr
&\stackrel{\iota}{\lra}&
\barr{c}
\bpi(10,30)(0,0)
\Gab \Gac \Gad \Gbc \Gbd
\put(10,10){\makebox(0,0){$\sqr$}}
\Gabbb \Gbbbc \Gbcbd
\Gacad
\epi
\earr
&\lra&
\raisebox{10pt}{
$\left(\barr{ccc}
\bpi(10,20)(0,0)
\Gaa \Gab \Gac \Gba \Gbb \Gbc
\Gaaba \Gbabb \Gbbbc
\Gabac
\epi
&\raisebox{35pt}{,}&
\bpi(10,20)(0,0)
\put(0,20){\makebox(0,0){$1$}}
\put(0,10){\makebox(0,0){$2$}}
\put(10,20){\makebox(0,0){$3$}}
\put(0,0){\makebox(0,0){$4$}}
\put(10,10){\makebox(0,0){$5$}}
\put(10,0){\makebox(0,0){$6$}}
\epi
\earr\right)$}
\earr
$$
\medskip
 
To define $\iota$ precisely, we will need to enlarge the set of tableaux
under consideration.  An {\it overlapping\/} rooted special rim-hook tableau
is a way of writing a shape as a union of special rim-hooks all of which are
disjoint with the exception of two whose intersection is precisely the
root, and this root must be a permissible element of both hooks.
Furthermore, one of the two hooks containing the root will be
designated as  {\it active}.  In our diagrams, the active hook will
have all of its nodes enlarged.  In a non-overlapping rooted special
rim-hook tableau, the unique hook containing the root is considered
active.

We are now ready to fully describe the involution $\iota$.  Starting with
an appropriately rooted special rim-hook tableau $\Sd$, we will
generate a sequence of tableaux
\beq
\label{seq}
\Sd=\Sd_0,\Sd_1,\Sd_2,\ldots,\Sd_t=\iota(\Sd).  
\eeq
For $0<i<t$, $\Sd_i$ will be an overlapping rooted special rim-hook tableau
of shape $\sh(\Sd)-r$ where $r$ is the root of $\Sd$.  Finally,
$\Sd_t$ will be a non-overlapping rooted special rim-hook tableau
satisfying~\ree{sh}.  All tableaux in the sequence will have the same type.
There are five operations for obtaining $\Sd_i$ from $\Sd_{i-1}$, 
$1\le i\le t$, as follows.  Whichever operation is applied, the
designation of ``active'' is transfered from the hook that was
modified to the new hook which it overlaps (except in the last step
where no new hook is overlapped and so the hook which was active
remains so).   Throughout, $r$ is the current root.
\ben
  \item[CO]  If $r$ is an internal or external corner of the active hook,
  reflect the root in a line containing the two nodes of the hook to which
  it is adjacent.
  \item[SI]  If the $r$ is in an active hook of size one, then it
  must be at the top or bottom of the first-column portion of the other hook
  containing it (since all roots are at permissible nodes).  Move the
  hook of size one and the root to the opposite end of that portion
  of the other hook.
  \item[HE] If $r$ is at the head of the active hook, then remove
  it and attach it just below the tail of this hook.
  \item[TV]  If $r$ is at the tail of the active hook of size at
  least two and the next node of the hook is directly above, then
  remove it and attach it to the head of the active hook.  Because of
  permissibility, exactly one of the two ways to attach the root
  (vertically or horizontally) will be possible.
  \item[TH]   If $r$ is at the tail of the active hook of size at
  least two and the next node of the hook is directly to the right,
  note that $r$ must also be at the tail of the other hook
  containing it (because of permissibility) and so the hooks must have
  different sizes (because they are special).  Let $s$ be the size of
  the smaller hook and let $v$ be the $(s+1)$st node from $r$ in the
  bigger hook.  Remove the portion of the bigger hook from $v$ on and
  attach it to the end of the smaller hook.  The root stays in the
  same place.
\een
Here is an example of the application of these rules to construct
$\iota(\Sd)$.  Each step is labeled with the operation being used.
When TH is applied, the vertex $v$ is marked.
\medskip
$$
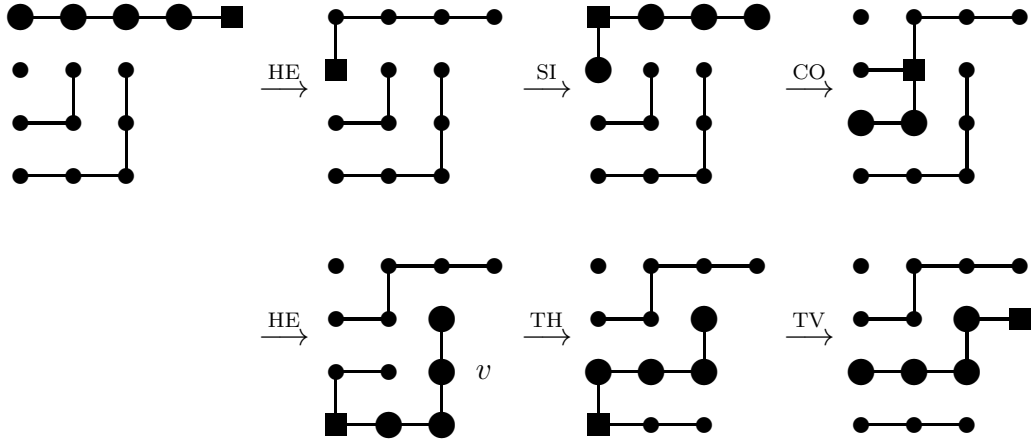
\thicklines
\setlength{\unitlength}{2pt}
\barr{ccccccc}
\bpi(40,30)(0,0)
\GadV{5} \GbdV{5} \GcdV{5} \GddV{5} \put(40,30){\makebox(0,0){$\sqr$}}
\Gac     \Gbc     \Gcc
\Gab     \Gbb     \Gcb
\Gaa     \Gba     \Gca
\Gadbd \Gbdcd \Gcddd \Gdded
\Gabbb \Gbbbc
\Gaaba \Gbaca \Gcacb \Gcbcc
\epi
&\raisebox{30pt}{$\stackrel{\rm HE}{\lra}$}&
\bpi(30,30)(0,0)
\Gad     \Gbd     \Gcd     \Gdd
\put(0,20){\makebox(0,0){$\sqr$}}
         \Gbc     \Gcc
\Gab     \Gbb     \Gcb
\Gaa     \Gba     \Gca
\Gadbd \Gbdcd \Gcddd \Gacad
\Gabbb \Gbbbc
\Gaaba \Gbaca \Gcacb \Gcbcc
\epi
&\raisebox{30pt}{$\stackrel{\rm SI}{\lra}$}&
\bpi(30,30)(0,0)
\put(0,30){\makebox(0,0){$\sqr$}}
         \GbdV{5} \GcdV{5} \GddV{5} 
\GacV{5} \Gbc     \Gcc
\Gab     \Gbb     \Gcb
\Gaa     \Gba     \Gca
\Gacad \Gadbd \Gbdcd \Gcddd 
\Gabbb \Gbbbc
\Gaaba \Gbaca \Gcacb \Gcbcc
\epi
&\raisebox{30pt}{$\stackrel{\rm CO}{\lra}$}&
\bpi(30,30)(0,0)
\Gad     \Gbd     \Gcd     \Gdd 
\Gac              \Gcc
         \put(10,20){\makebox(0,0){$\sqr$}}
\GabV{5} \GbbV{5} \Gcb
\Gaa     \Gba     \Gca
\Gacbc \Gbcbd \Gbdcd \Gcddd 
\Gabbb \Gbbbc
\Gaaba \Gbaca \Gcacb \Gcbcc
\epi
\\[30pt]
&\raisebox{30pt}{$\stackrel{\rm HE}{\lra}$}&
\bpi(30,30)(0,0)
\Gad     \Gbd     \Gcd     \Gdd 
\Gac     \Gbc     \GccV{5}
\Gab     \Gbb     \GcbV{5} 
\put(28,10){\makebox(0,0){$v$}}
\put(0,0){\makebox(0,0){$\sqr$}}
         \GbaV{5} \GcaV{5}
\Gacbc \Gbcbd \Gbdcd \Gcddd 
\Gaaab \Gabbb 
\Gaaba \Gbaca \Gcacb \Gcbcc
\epi
&\raisebox{30pt}{$\stackrel{\rm TH}{\lra}$}&
\bpi(30,30)(0,0)
\Gad     \Gbd     \Gcd     \Gdd 
\Gac     \Gbc     \GccV{5}
\GabV{5} \GbbV{5} \GcbV{5} 
\put(0,0){\makebox(0,0){$\sqr$}}
         \Gba     \Gca
\Gacbc \Gbcbd \Gbdcd \Gcddd 
\Gaaab \Gabbb \Gbbcb \Gcbcc
\Gaaba \Gbaca 
\epi
&\raisebox{30pt}{$\stackrel{\rm TV}{\lra}$}&
\bpi(30,30)(0,0)
\Gad     \Gbd     \Gcd     \Gdd 
\Gac     \Gbc     \GccV{5} \put(30,20){\makebox(0,0){$\sqr$}}
\GabV{5} \GbbV{5} \GcbV{5} 
\Gaa     \Gba     \Gca
\Gacbc \Gbcbd \Gbdcd \Gcddd 
\Gabbb \Gbbcb \Gcbcc \Gccdc
\Gaaba \Gbaca 
\epi
\earr
$$
\medskip

We must now show that $\iota$ is a well-defined sign-reversing
involution.  By examining each of the rules in turn it is easy to see
that all the rim-hook tableaux generated have the same type and that
the root is always in a permissible cell.  Also it is clear that each rule
can be reversed:  CO, SI, and TH are self-inversive while HE and TV are
inverses of each other.  So if we do a single step without changing
which hook is designated as active and then apply the rules again
we will return to the original tableau.  This means that $\iota$ will
be an involution if it is well defined.

To finish the proof of well definedness, we must show that the
algorithm terminates, i.e., that eventually a non-overlapping tableau
is produced.  Suppose to the contrary that~\ree{seq} goes on forever.
Since there are only a finite number of tableaux of a given type, this
sequence must repeat.  Let $\Sd_i$ have the smallest index such that
$\Sd_i=\Sd_j$ for some $j>i$.  Note that since we are assuming that
all tableaux other than $\Sd_0$ are overlapping, we must have $i\ge1$.
But we saw in the previous paragraph that each step of the algorithm
is invertible and $i\ge1$, so we must have $\Sd_{i-1}=\Sd_{j-1}$.  This
contradicts the minimality of $i$ and proves termination.

It remains to show that $\iota$ is sign reversing.  In fact, we will
prove the stronger statement that the sign of $\Sd_i$ depends only on
the type of active hook it contains and the sign of $\Sd$. 
We will label the active hook in a tableau with the same pair of
letters used for the rule which can be applied to it. 
We will also have to refine CO into
two types of hooks, namely CI (respectively, CE) if the corner is internal
(respectively, external).  So in the previous example, the CO hook in the
third tableau would be classified as a CI.
Similarly, we have to split HE into
two types:  HV (respectively, HH) if the head of the hook
is directly above (respectively, to the right of) its predecessor in
the hook.  Looking at the previous example again, the initial active
hook is an HH while the fourth one is an HV.
\ble
Let $H_i$ be the active hook in $\Sd_i$ from the sequence~\ree{seq}
where $0\le i<t$.  Then
$$
\ep(\Sd_i)=\case{\ep(\Sd)}{if $H_i$  is {\rm CE}, {\rm HH}, or {\rm TV},}
{-\ep(\Sd)}{if $H_i$ is {\rm CI}, {\rm HV}, or {\rm TH}.}
$$
\ele   
\proof
We induct on $i$.  Suppose first that $i=0$.  Since $n$ must be at
the end of a row and column, the only possibilities for $H_0$ are CE,
HH, or TV.  Clearly we also have $\ep(\Sd_0)=\ep(\Sd)$, and so the
lemma is true in this case.

The induction step breaks down into six cases depending upon the
nature of $H_i$.  Since they are all similar, we will do the one for
$H_i$ = TV and then the other five can be verified by the reader if
they wish.  

Note first that because of the way the root moves in each of the
steps, TV can only be preceded by CE, TH, TV, or SI.  (For example,
after a CI step the root cannot be in the first column and so it can't
precede TV.)  

If $H_{i-1}$ is CE then by induction and the fact that 
a CE step does not change the number of vertical edges in a hook 
$$
\ep(\Sd_i)=\ep(\Sd_{i-1})=\ep(\Sd).
$$

If $H_{i-1}$ is TH then this step changes the number of vertical
edges in the overlapping hooks by $\pm1$.  So by induction again
$$
\ep(\Sd_i)=-\ep(\Sd_{i-1})=-(-\ep(\Sd))=\ep(\Sd).
$$

If $H_{i-1}$ is TV, then in order for it to overlap $H_i$ by a node in
the first column it must be that $H_{i-1}$ consists precisely of a strip
of nodes in the first column directly below the root of $H_i$.  So the
TV step does not change the number of vertical edges in $H_{i-1}$ and
so the same equalities as for the CE case give the desired conclusion.

Finally, consider what happens if $H_{i-1}$ is SI.   Then in
$\Sd_{i-1}$, the singleton hook must be at the top of the portion
of $H_i$ in the first column.  So $H_{i-2}$ must have been either a CE
(if $H_i$ contains nodes outside of the first column) or a TV of the type
discussed in the previous paragraph.  In either case, the passage from
$H_{i-2}$ to $H_{i-1}$ causes no change in sign and neither does the
SI step.  So by induction
$$
\ep(\Sd_i)=\ep(\Sd_{i-2})=\ep(\Sd).
$$

This completes the demonstration that the lemma holds in the TV
case. \qed

The follow theorem will complete the proof that $\iota$, and hence
$I$, is a well-defined, sign-reversing involution.
\bth
If $\iota(\Sd)=\Sd'$ then $\ep(\Sd')=-\ep(\Sd)$
\eth
\proof
We continue to use the notation in~\ree{seq} where $\Sd'=\Sd_t$.
Since $\Sd_t$ is non-overlapping, $H_t$ must be one of
CE, HH, or TV.  

If $H_t$ is CE then $H_{t-1}$ must have been CI.  So using the
previous lemma with $i=t-1$ and the fact that a CI step does not
change sign gives
$$
\ep(\Sd')=\ep(\Sd_{t-1})=-\ep(\Sd).
$$

If $H_t$ is HH then $H_{t-1}$ must have been TV.  Furthermore, this TV
step must have added a horizontal step to $H_{t-1}$ and so decreased
the number of vertical edges by one.  So using the lemma again gives 
$$
\ep(\Sd')=-\ep(\Sd_{t-1})=-\ep(\Sd).
$$

Finally, suppose $H_t$ is TV.  Then $H_{t-1}$ is either HH or HV and
either the previous string of equalities or the one before
that, respectively, hold.  So in all cases $\ep(\Sd')=-\ep(\Sd)$. \qed

\section{The \tpo-free Conjecture}
\label{tfc}

In order to make a connection of our work with the  \tpo-free
Conjecture, we first need to introduce Stanley's chromatic symmetric
function~\cite{sta:sfg,sta:gcr}.  Let $G=(V,E)$ be a graph with
a finite set of vertices $V$ and edges $E$.  A {\it proper coloring\/}
of $G$ from a set $A$ is a function $\ka:V\ra A$ such that $uv\in E$ implies
$\ka(u)\neq\ka(v)$.  Now consider a countably infinite set of
variables $\bx=\{x_1,x_2,\ldots\}$.  Stanley associated with each
graph a formal power series
$$
X_G= 
X_G(\bx)=\sum_{\kappa:V\ra\bbP} x_{\ka(v_1)}x_{\ka(v_2)}\cdots x_{\ka(v_n)}
$$
where $\ka$ is a proper coloring from the positive integers $\bbP$.
Note that if one sets $x_1=x_2=\ldots=x_n=1$ and $x_i=0$ for $i>n$,
denoted $\bx=1^n$,
then $X_G$ reduces to the number of proper colorings of $G$ from a set
with $n$ elements.  So under this substitution, $X_G(1^n)=P_G(n)$
where $P_G(n)$ is the famous chromatic polynomial of
Whitney~\cite{whi:lem}.   Also, because permuting the colors of a
proper coloring keeps the coloring proper, $X_G(\bx)$ is in the
algebra $\La(\bx)$ of symmetric functions in $\bx$ over the rationals.
In~\cite{sta:sfg,sta:gcr},
Stanley was able to derive many interesting properties of the
chromatic symmetric function $X_G(\bx)$ some of which generalize those
of the chromatic polynomial and some of which cannot be interpreted
after substitution.  

One natural question to ask is whether one can say anything about the
expansion of $X_G(\bx)$ in any of the usual bases for $\La(\bx)$.  If
$f\in\La(\bx)$ and $\{b_\la\ :\ \la\in\Par\}$ is a basis for
$\La(\bx)$, then we will say that $f$ is {\it $b$-positive\/} if in the
expansion $f=\sum_\la c_\la b_\la$, all of the coefficients satisfy
$c_\la\ge0$.  The \tpo-free Conjecture states that for certain graphs,
$X_G$ is $e$-positive where the $e_\la$ are the elementary symmetric
functions.

To describe the appropriate graphs for the conjecture, consider a
finite poset (partially ordered set) $(P,\le)$.  We say that $P$ is 
{\it ${\rm (}{\bf a}{\rm +}{\bf b}{\rm )}$-free\/} if it contains no
induced subposet isomorphic to a disjoint union of an $a$-element
chain and a $b$-element chain.  Also, given any poset $P$, we can form its
{\it incomparability graph}, $G(P)$, having vertices $V=P$ and an edge
between $u$ and $v$ in $G(P)$ if and only if $u$ and $v$ are
incomparable in $P$.   Through their work on immanants of Jacobi-Trudi
matrices, Stanley and Stembridge~\cite{ss:ijt} were led  to the
following conjecture.
\bcon[\tpo-free Conjecture]
Let $P$ be a a \tpo-free poset.  Then $X_{G(P)}$ is $e$-positive,
i.e., if
\beq
\label{Xe}
X_{G(P)}=\sum_\mu c_\mu e_\mu
\eeq
then $c_\mu\ge0$ for all $\mu$.
\econ 

There is a fair amount of evidence  to support this conjecture.
Stembridge has verified that it
is true for all 884 \tpo-free posets having at most 7 elements.
Gebhard and Sagan~\cite{gs:csf} have used the theory of symmetric
functions in noncommuting variables to prove that the conjecture holds
for certain posets which are both \tpo- and 
${\rm (}{\bf 2}{\rm +}{\bf 2}{\rm )}$-free.  

One of the most significant results about this conjecture was obtained
by Gasharov~\cite{gas:igt} who proved that if $P$ is \tpo-free then
$X_{G(P)}$ is $s$-positive where $s_\la$ is the Schur function
corresponding to $\la$.  He did this by giving a combinatorial
interpretation to the coefficients in the $s$-expansion of $X_{G(P)}$
which we will need in the sequel.

For a poset $P$, a  {\it $P$-tableau $T$ of
shape $\la$\/} is a filling of the cells of $\la$ with the elements of
$P$ (each used exactly once) such that for all $(i,j)\in\la$:
\ben
\item $T_{i,j}<T_{i+1,j}$, and
\item $T_{i,j}\not> T_{i,j+1}$
\een
where a condition is considered vacuously true if subscripts refer to
a cell outside of $\la$.  Note that when $P$ is a chain, then a
$P$-tableau is just a standard Young tableau.  Letting $f_P^\la$
denote the number of $P$-tableaux of shape $\la$, Gasharov proved the
following result.
\bth[Gasharov]
If $P$ is \tpo-free then
\beq
\label{Xs}
X_{G(P)}=\sum_\la f_P^{\la} s_{\la'}
\eeq
where $\la'$ is the conjugate of $\la$. \qed
\eth
Note that this immediately implies $s$-positivity.  

Chow~\cite{cho:nci} pointed out that~\ree{Xs} could be combined with 
E\~gecio\~glu and Remmel's result to obtain a
combinatorial interpretation of the coefficients $c_\mu$ in~\ree{Xe}. 
First note that the change of basis matrix between the Schur and
elementary symmetric functions is  
$$
s_{\la'}=\sum_\mu K^{-1}_{\mu,\la} e_\mu
$$
Combining this with~\ree{Xs} we get
$$
X_{G(P)}=\sum_{\la,\mu} K^{-1}_{\mu,\la} f_P^{\la} e_\mu.
$$
Since the $e_\mu$ are a basis, we have
$$
c_\mu=\sum_{\la} K^{-1}_{\mu,\la} f_P^{\la}.
$$
Finally we apply Theorem~\ref{K-1} to get the desired interpretation.
\bco[Chow]
\label{cmu}
The coefficients $c_\mu$ in the $e$-expansion of $X_G(P)$ satisfy
$$
c_\mu=\sum_{(S,T)} \ep(S)
$$
where the sum is over all pairs of a special rim hook tableau $S$ of
type $\mu$ and a $P$-tableau $T$ with the same shape as $S$.  \qed
\eco
Note that a column of a $P$-tableau $T$ must be a chain in
$P$ and the number of rim hooks in $S$ is at most the length of its
first column because they are special.  So
the previous corollary implies that $c_\mu=0$ whenever $\mu$ has more
parts than the height of $P$, $h(P)$ (which is defined as the number
of elements in the longest chain of $P$).  So this cuts down on the
number of coefficients which we need to consider.

So to show that $c_\mu\ge0$, it suffices to find an involution $I$ on
the  pairs in the previous corollary such that if $I(S,T)=(S',T')$ then
\ben
\item $(S,T)=(S',T')$ implies $\ep(S)=1$, and
\item $(S,T)\neq(S',T')$ implies $\ep(S')=-\ep(S)$.
\een
The involution in the  previous section does this when
$P$ is a chain (so that $P$-tableaux are standard Young tableaux).
Of course, then $G(P)$ is a totally disconnected graph and
thus $X_{G(P)}=e_{(1^n)}$ directly from the definition of $X_G$.  So it
would be nice to apply these ideas to a less obvious case.  This will
be done in the next theorem which was also derived by Stanley and
Stembridge themselves~\cite{ss:ijt} using the theory of rook
placements.  
Note that if $h(P)\le2$ then $P$ must be \tpo-free since it does not
even contain a 3-element chain.

\bth[Stanley-Stembridge]
If $h(P)\le2$ then $X_{G(P)}$ is $e$-positive.
\eth
\proof 
It suffices to construct the involution $I$ in this case.  First of
all, the remarks after Corollary~\ref{cmu} show 
that any $P$-tableau, $T$, will have at most two rows
since $h(P)\le2$.  If a
special rim-hook tableau $S$ has at most two rows and $\ep(S)=-1$,
then it must have exactly two rows and contain exactly one vertical
edge.  

So let $(S,T)$ be a pair having common shape $\la=(\la_1,\la_2)$ and
with $\ep(S)=-1$.  Let $\Sd$ be the tableau obtained
by rooting $S$ at the end of the second row.  Now form
$\Sd'=\iota(\Sd)$.  By our assumptions on $S$, $\iota$ will consist of
a sequence of steps resulting in an $\Sd'$ which will be the unique tableau of
shape $(\la_1+1,\la_2-1)$ where each row is a special rim hook and the
root is at the end of the first row.  We now obtain $(S',T')=I(S,T)$
as before by removing the root from $\Sd'$ to get $S'$ and getting
$T'$ by leaving all the elements of $T$ in the same place except the
one, $x$, which was at the end of the second row in $T$ and is now at
the end of the  first row in $T'$.  Clearly this reverses sign, but we
must check that $T'$ is still a valid $P$-tableau.  But $x$ was at the
end of a column of length two in $T$ and so cannot have any element
above it in $P$ which is of height two.  So the row condition for
$P$-tableaux will still be satisfied and $I$ is well defined. 

In order to turn $I$ into an involution, we must characterize those
pairs $(S',T')$ which are in the image of the function so far and map
them back to their preimages.  So suppose the tableaux in $(S',T')$
have shape $\nu=(\nu_1,\nu_2)$ where $\nu_1>\nu_2+1>0$ and
$\ep(S')=+1$.   Consider the elements $x=T'_{1,\nu_1}$ and
$y=T'_{1,\nu_2+1}$.  Then it is easy to see that the image of $I$ is
precisely the set of all such pairs such that $x>y$ in $P$.  So we can
reverse the algorithm by forming $\Sd'$ which is $\Sd$ rooted at the
end of the first row, applying $\iota$, and then using the usual
operations to recover $(S,T)$.

Finally, we need to consider what happens to the $(S,T)$ which
have not been paired up so far by $I$.
But since these all have positive
sign, we can just make them fixed points of $I$.  This completes the
definition of $I$ and the proof of the Theorem.
\qed
\medskip

As an example of the algorithm in the previous proof, consider the
poset
\medskip
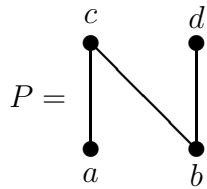
\bce
\thicklines
\setlength{\unitlength}{2pt}
\bpi(30,40)(-10,-10)
\put(-10,10){\makebox(0,0){$P=$}}
\Gaa \put(0,-5){\makebox(0,0){$a$}}
\Gca \put(20,-5){\makebox(0,0){$b$}}
\Gac \put(0,25){\makebox(0,0){$c$}}
\Gcc \put(20,25){\makebox(0,0){$d$}}
\Gaaac \Gcaac \Gcacc
\epi
\ece
\medskip
There are 20 pairs $(S,T)$ for this poset.  
To present them in an economical way, we will combine each pair into a
single tableau with elements in the same places as in $T$ and edges
between pairs of elements which are adjacent in a hook of $S$.
With this notation, tableaux which are matched by $I$ are as follows.
$$
\barr{ccccccccccccc}
a&\h&b&\h&d&\h&c&\llra&a&\h&b&\h&d\\
 &  & &  & &  & &     &|\\
 &  & &  & &  & &     &c\\[10pt]
b&\h&a&\h&c&\h&d&\llra&b&\h&a&\h&c\\
 &  & &  & &  & &     &|\\
 &  & &  & &  & &     &d\\[10pt]  
b&\h&a&\h&d&\h&c&\llra&b&\h&a&\h&d\\
 &  & &  & &  & &     &|\\
 &  & &  & &  & &     &c\\[10pt]  
b&\h&d&\h&a&\h&c&\llra&b&\h&d&\h&a\\
 &  & &  & &  & &     &|\\
 &  & &  & &  & &     &c\\[10pt]  
a&\h&b&\h&d&  & &\llra&a&  &b\\
 &  & &  & &  & &     & &  &|\\
c&  & &  & &  & &     &c&\h&d\\[10pt]  
b&\h&a&\h&c&  & &\llra&b&  &a\\
 &  & &  & &  & &     & &  &|\\
d&  & &  & &  & &     &d&\h&c
\earr
$$
And here is a list of the fixed points organized in columns by shape.
$$
\barr{ccccccccccccccccc}
a&\h&b&\h&c&\h&d&\rule{30pt}{0pt}&b&\h&a&\h&d&\rule{30pt}{0pt}&a&\h&b\\
 &  & &  & &  & &                & &  & &  & &                & &  &\\ 
 &  & &  & &  & &                &c&  & &  & &                &c&\h&d\\[15pt] 
b&\h&c&\h&d&\h&a&\rule{30pt}{0pt}&b&\h&d&\h&a&\rule{30pt}{0pt}&b&\h&a\\
 &  & &  & &  & &                & &  & &  & &                & &  &\\ 
 &  & &  & &  & &                &c&  & &  & &                &d&\h&c\\[15pt] 
c&\h&d&\h&a&\h&b&                & &  & &  & &                & &  &\\[20pt] 
d&\h&a&\h&b&\h&c&                & &  & &  & &                & &  &
\earr
$$
Counting the fixed points by shape, we immediately have
$$
X_{G(P)}=4e_{(4)}+2e_{(3,1)}+2e_{(2,2)}.
$$

We remark that if one adjoins a unique maximum or minimum to $P$, then
this adds an isolated vertex to $G(P)$.  So this just multiplies
$X_{G(P)}$ by $e_{(1)}$.  So the previous theorem implies that $X_{G(P)}$
is $e$-positive in the case when $h(P)=3$ and $P$ has a unique maximum
or minimum as well as in the case when $h(P)=4$ and $P$ has both.

\section{Further work}
\label{fw} 

We hope that the involution we have presented will just be a first
step towards making progress on the problem of E\~gecio\~glu and
Remmel as well as on the \tpo-free Conjecture.  In order to encourage
the reader to develop these ideas, let us present some thoughts about
how to proceed.

To complete a combinatorial proof of  $K^{-1}K=I$, we must find an
involution $I$ on pairs $(S,T)$ of the same shape for any given type
$t(S)=\mu$ and content $c(T)=\nu$.  Of course, $I$ should be sign
reversing on all its $2$-cycles and any fixed points should have $S$ of
positive sign.  If $\mu\neq\nu$ then there should be no fixed points.
If $\mu=\nu$ then there should be a single fixed point which should
probably be the unique pair of shape $\mu=\nu$.

If the largest element of $T$ to appear in a hook of size at least two
only occurs once in $T$, then one can use the same algorithm as before
to construct $I$.  The question is what to do if the largest element
of $T$ occurs with nontrivial multiplicity.  One possible solution is
to recall that to every semistandard $T$ there is a canonically
associated standard Young tableau $T_0$ obtained by labeling the $1$'s
in $T$ from left to right with $1,2,\ldots,\nu_1$, then labeling the
$2$'s in the same manner with $\nu_1+1,\nu_1+2,\ldots,\nu_1+\nu_2$,
and so forth.  One can now apply the old involution $I$ to the pair
$(S,T_0)$ to obtain a pair $(S',T_0')$.  The problem is that if we now
reverse the standardization procedure, we may no longer get a
semistandard tableau as one may get two largest elements in the same
column.  But perhaps there is a way to circumvent this.

Another possible approach would be to try and come up with a new set
of rules where all the largest elements were moved as a group in each
step.  We have been able to see how some pairs might behave under
this assumption, but have not been able to come up with something that
works in all circumstances.

The difficulties when dealing with posets are similar.  Since $P$ may
have many maximal elements, it is unclear which of them should be used
to root $S$.  Or maybe it is the case that the root can be chosen
arbitrarily among these maximals as long as one always chooses the
same element.  But some extra idea  will have to be incorporated to
deal with the fact that a maximal element need not be at the end of a
row (although it must be at the end of a column) and to
make sure that when the maximal element is moved to the new position,
the new array remains a $P$-tableau.  This should be where the
\tpo-free condition comes in.

\begin{\bib}{99}

\bibitem{cho:nci} T. Chow, A note on a combinatorial interpretation of
the  $e$-coefficients of the chromatic symmetric function, preprint
(1997), 9 pp.

\bibitem{er:cii} \"O. E\~gecio\~glu and J. Remmel, A combinatorial
interpretation of the inverse Kostka matrix, {\it Linear and Multilinear
Algebra\/} {\bf 26} (1990), 59--84.

\bibitem{ful:yt} W. Fulton, ``Young Tableaux,'' London Mathematical
Society Student Texts 35, Cambridge University Press, Cambridge, 1999.

\bibitem{gs:csf} D. Gebhard and B. Sagan, A chromatic symmetric
function in noncommuting variables, {\it J. Algebraic Combin.} {\bf
13} (2001), 227--255. 

\bibitem{gas:igt} V. Gasharov, Incomparability graphs of ({\bf 3}+{\bf
1})-free posets are  $s$-positive, {\it Discrete Math.} {\bf 157}
(1996), 193--197.

\bibitem{gv:bdp} I.  Gessel  and  G.  Viennot,  Binomial 
determinants, paths, and hook length formulae, {\it \aim} {\bf 58}
(1985), 300--321.

\bibitem{gv:dpp} I. Gessel and G. Viennot, Determinants, paths, and
plane partitions, in preparation.

\bibitem{lin:vri} B. Lindstr\"om, On the vector representation of induced
matroids, {\it Bull.\  London Math.\  Soc.\ } {\bf 5} (1973), 85--90.

\bibitem{mac:sfh} I. G. Macdonald,  ``Symmetric  functions 
and Hall polynomials,'' 2nd edition, \oup, Oxford, 1995.

\bibitem{sag:sg} B. Sagan, ``The Symmetric Group: Representations,
Combinatorial Algorithms,  and Symmetric Functions,'' 2nd edition,
Springer-Verlag, New York, 2001. 

\bibitem{sta:sfg} R. P. Stanley, A symmetric function generalization
of the chromatic polynomial of a graph, {\it Advances in Math.}\ 
{\bf 111} (1995), 166--194.

\bibitem{sta:gcr} R. P. Stanley, {\bf G}raph Colorings {\bf a}nd
{\bf r}elated {\bf s}ymmetric functions: {\bf i}deas and 
{\bf a}pplications: A description of results, interesting applications,
\& notable open problems,  Selected papers in honor of Adriano Garsia
(Taormina, 1994) {\it Discrete Math.\/} {\bf 193} (1998), 267--286.

\bibitem{sta:ec2} R. P. Stanley, ``Enumerative Combinatorics,
Volume 2,''  Cambridge University Press, Cambridge, 1999.

\bibitem{ss:ijt} R. P. Stanley and J. Stembridge, On immanants of
Jacobi-Trudi  matrices and permutations with restricted position, 
{\it J. Combin.\ Theory Ser.\ A}  {\bf 62} (1993), 261--279.

\bibitem{sw:scr} D. Stanton and D. White, A Schensted correspondence for  rim 
hook tableaux, {\it \jcta} {\bf 40} (1985), 211--247.

\bibitem{whi:bpo} D.  White,   A  bijection  proving  orthogonality  of  the 
characters of $S_n$,  {\it \aim} {\bf 50} (1983), 160--186.

\bibitem{whi:lem}  H. Whitney, A logical expansion in mathematics,
{\it Bull.\ Amer.\ Math.\ Soc.\ } {\bf 38} (1932), 572--579.

\end{\bib}

\end{document}